\documentclass{amsart}

\theoremstyle{plain}
\newtheorem{thm}{\indent\sc Theorem}[section]
\newtheorem{prop}[thm]{\indent\sc Proposition}

\theoremstyle{definition}
\newtheorem{defn}[thm]{\indent\sc Definition}

\theoremstyle{remark}

\def\proofwp{{\sc Proof}}

\numberwithin{equation}{section}
\def\cA{\mathcal{A}}

\def\cC{\mathcal{C}}

\usepackage{amscd}
\usepackage{amsmath,amssymb}
\usepackage{amssymb}
\usepackage{amsmath}
\usepackage{amsthm}
\usepackage{latexsym}
\usepackage{amsfonts}
\usepackage{color}
\usepackage{graphics}
\usepackage[all,knot]{xy}
\usepackage[all]{xy}
\xyoption{arc}
\usepackage{enumerate}
\usepackage{amstext}

\def\s{{\bf s}}

\def\im{{\hbox{\bf Im}\;}}

\def\bd{\mbox{\bf bd }}

\def\mapr#1{\smash{\mathop{\buildrel{#1}\over\longrightarrow}}}

\newcommand{\RR}{\mathbb{R}}

\begin{document}

\title[Mayer-Vietoris sequence in cohomology of Lie algebroids]{Mayer-Vietoris sequence in cohomology of Lie algebroids on simplicial complexes}



\author{{\textsc{Jose R. Oliveira}}}
\address{Department of Mathematics, Minho University,
    Braga, Portugal}
\curraddr{}
\email{jmo@math.uminho.pt}
\thanks{The author was partially supported by MICINN, Grant MTM2014-56950-P}



\date{}

\dedicatory{}


\maketitle






\begin{center}

\vspace{3mm}

\textsc{Abstract}

\vspace{3mm}

\end{center}


It is shown that the Mayer-Vietoris sequence holds for the cohomology of complexes of Lie algebroids which are defined on simplicial complexes and satisfy the compatibility condition concerning restrictions to the faces of each simplex. The Mayer-Vietoris sequence will be obtained as a consequence of the extension lemma for piecewise smooth forms defined on complexes of Lie algebroids.

\vspace{3mm}







\section{Introduction}

\vspace{3mm}

D. Sullivan in \cite{suli-inf} and H. Whitney in \cite{wity-git} considered several cell-like constructions of cochain complexes which induce isomorphisms in cohomology with classical cohomologies of the underlying polytope. Mishchenko and Oliveira in \cite{mish-oli} have extended Sullivan-Whitney constructions to transitive Lie algebroids defined over triangulated manifolds. The main construction considered by Mishchenko and Oliveira in \cite{mish-oli} is the one in which a transitive Lie algebroid on a triangulated smooth manifold is fixed and one considers the family of Lie algebroids obtained by restriction of the Lie algebroid to each simplex of the base. The notion of piecewise smooth form is defined in a similar way to Whitney forms on a cell space. In \cite{mish-oli}, it is proved that the cohomology of this construction is isomorphic to the cohomology of the Lie algebroid. The Mayer-Vietoris sequence corresponding to two open subsets of the manifold obtained by the union of open stars is used in the course of the proof of that result.

The work synthesized in the present paper is the extension of the Mayer-Vietoris property for the cohomology of families of transitive Lie algebroids defined over simplicial complexes. For this purpose, we follow the main construction presented by Mishchenko and Oliveira in their paper \cite{mish-oli}. Based in that construction, we are going to consider, in this paper, families of transitive Lie algebroids defined over simplicial complexes, satisfying a compatibility property concerning restrictions of the Lie algebroids to the faces of each simplex. When this structure is given, we will define piecewise smooth, following Mishchenko-Oliveira's work (\cite{mish-oli}). Furthermore, a differential can be defined, yielding a commutative differential graded algebra. Its cohomology is, by definition, the piecewise smooth cohomology of the family of Lie algebroids considered. The aim of the present paper is to state and prove the Mayer-Vietoris sequence for the cohomology corresponding to this construction. We will show that the extension lemma for piecewise smooth forms on such families of transitive Lie algebroids is an essential key on the proof of short exact sequence which generates the Mayer-Vietoris sequence in cohomology.

Throughout the paper, all manifolds considered are smooth, finite-dimensional and possibly with boundaries of different indices.


\vspace{3mm}

\textbf{Acknowledgments}. The author wishes to thank to Aleksandr Mishchenko, Jesus Alvarez, Nicolae Teleman and James Stasheff for their strong dynamism to discuss several topics concerning cohomology of cell spaces and Lie algebroids and also to the referee for many helpful comments and suggestions.


\vspace{6mm}

\begin{center}
\section{Complexes of Lie algebroids and differential forms}
\end{center}

\vspace{3mm}

We deal with transitive Lie algebroids defined over the simplices of a simplicial complex. We briefly discus a class of spaces for which piecewise smooth cohomology is defined, precisely the class of all complexes of Lie algebroids. In what follows, all simplicial complexes considered are geometric and finite. Simplex means always closed simplex and each simplex can be represented as a convex body generated by its vertices. We shall denote the boundary of the simplex $\Delta$ by $\bd\Delta$. We shall write $\s\prec\Delta$ to indicate that $\s$ is a face of the simplex $\Delta$. The notation $\varphi:\s\hookrightarrow \Delta$, in which $\varphi$ is the inclusion, will also be used when $\s$ is a face of $\Delta$. Various definitions and properties stated through the entire paper can be found, on level of cell spaces, in \cite{jesus-cal}, \cite{bott-dft}, \cite{Eile-stin}, \cite{suli-inf}, and \cite{wity-git}.

We begin by reviewing basic definitions and constructions concerning Lie algebroids. An extensive discussion on these issues can be found in the book \cite{makz-lga} by Mackenzie. The papers \cite{mish-oli}, \cite{mish-cla} and \cite{mish-latb} contain a summary of those issues.


\vspace{3mm}

Let $M$ be a smooth manifold, possibly with boundaries of different indices, $TM$ the tangent bundle to $M$ and $\Gamma(TM)$ the Lie algebra of the vector fields on $M$. We recall that a Lie algebroid on $M$ (\cite{makz-lga}, \cite{mish-oli}) is a vector bundle $\pi:\cA\longrightarrow M$ on $M$ equipped with a vector bundle morphism $\gamma:\cA\longrightarrow TM$, called anchor of $\cA$, and a structure of real Lie algebra on the vector space $\Gamma(\cA)$ of the sections of $\cA$ such that the map $\gamma_{\Gamma}:\Gamma(\cA)\longrightarrow\Gamma(TM)$, induced by $\gamma$, is a Lie algebra homomorphism and the action of the algebra $\cC^{\infty}(M)$ on $\Gamma(\cA)$ satisfies the natural condition:
$$[\xi,f\eta]=f[\xi,\eta] + (\gamma_{\Gamma}(\xi)\cdot f)\eta$$ for each $\xi$, $\eta$ $\in \Gamma(\cA)$ and $f\in \cC^{\infty}(M)$. The Lie algebroid $\cA$ is called transitive if the anchor $\gamma$ is surjective. We recall that Lie algebroids are the infinitesimal objects of Lie groupoids (for the definition of Lie groupoid, see Mackenzie \cite{makz-lga}, definition 1.1.1 and definition 1.1.3). Lie groupoids and Lie algebroids enjoy some of the properties of Lie groups and Lie algebras. We notice that not every Lie algebroid is integrable to a Lie groupoid. The theorem 4.1 of the paper \cite{crac-loja} shows necessary and sufficient conditions so that a Lie algebroid is integrable to a Lie groupoid.

\vspace{3mm}

Let $M$ be a smooth manifold and suppose that $\cA$ is a transitive Lie algebroid on $M$. Let $\varphi:N\hookrightarrow M$ be a submanifold, possibly with boundaries of different indices. We recall that the Lie algebroid restriction of $\cA$ to the submanifold $N$, denoted by $\cA^{!!}_{N}$, is the Lie algebroid $\varphi^{!!}\cA$ constructed as inverse image of $\cA$ by the mapping $\varphi$ (details of image inverse of Lie algebroids can be seen in the section 4.2 of \cite{makz-lga} or in the second section of \cite{mish-oli}). The Lie algebroid $\cA^{!!}_{N}$ is transitive.


\vspace{3mm}

\begin{defn} Let $K$ be a simplicial complex. A \textit{complex of Lie algebroids} on $K$ is a family $\underline{\cA}=\{\cA_{\Delta}\}_{\Delta\in K}$ such that, for each $\Delta\in K$, $\cA_{\Delta}$ is a transitive Lie algebroid on $\Delta$ and $\cA_{\Delta'}=(\cA_{\Delta})^{!!}_{\Delta'}$ for each face $\Delta'$ of $\Delta$, that is, the Lie algebroid restriction of $\cA_{\Delta}$ to $\Delta'$ is the Lie algebroid $\cA_{\Delta'}$.
\end{defn}


Alternatively, a complex of Lie algebroids on $K$ means a family of transitive Lie algebroids defined on the simplices of $K$ such that the structures of Lie algebroids induced on each intersection of two simplices coincide.


\begin{defn} Let $K$ be a simplicial complex and $\underline{\cA}=\{\cA_{\Delta}\}_{\Delta\in K}$ a complex of Lie algebroids on $K$. Let $L$ be a simplicial subcomplex of $K$. We can consider a new complex of Lie algebroids, defined on $L$ and denoted by $\underline{\cA}^{L}$, given by restriction of $\underline{\cA}$ to the simplices of $L$, that is, $\underline{\cA}^{L}=\{\cA_{\Delta}\}_{\Delta\in L}$. The complex $\underline{\cA}^{L}$ is called the restriction of the complex of Lie algebroids $\underline{\cA}$ to the simplicial subcomplex $L$.
\end{defn}


We give now some examples of complexes of Lie algebroids.

\vspace{3mm}

\textit{Example 1} (Tangent complex). Let $K$ be a simplicial complex. For each simplex $\Delta\in K$, consider the tangent Lie algebroid $T\Delta$ defined over $\Delta$. If $\Delta'$ is a face of $\Delta$ then $(T\Delta)^{!!}_{\Delta'}=T\Delta'$ (see proposition 2.6 of \cite{mish-oli}) and consequently we obtain a complex of Lie algebroids $\{T\Delta\}_{\Delta\in K}$, which is called the corresponding tangent complex on $K$.

\vspace{3mm}

\textit{Example 2} (Trivial complex). Let $K$ be a simplicial complex and $\frak g$ a real finite dimensional Lie algebra. For each simplex $\Delta\in K$ consider the transitive Lie algebroid $T\Delta \oplus (\Delta\times \frak g)$. If $\Delta'$ is a face of $\Delta$ then
$(T\Delta \oplus (\Delta\times \frak g))^{!!}_{\Delta'}=T\Delta'\oplus (\Delta'\times \frak g)$ (see proposition 2.7 of \cite{mish-oli}). We conclude that the family $\{T\Delta \oplus (\Delta\times \frak g)\}_{\Delta\in K}$ is a complex of Lie algebroids. This complex is called the trivial complex on $K$.

\vspace{3mm}

\textit{Example 3} (Complex corresponding to a triangulated manifold). Let $M$ be a smooth manifold, smoothly triangulated by a simplicial complex $K$, and $\cA$ a transitive Lie algebroid on $M$. We can consider the Lie algebroid restriction $\cA^{!!}_{\Delta}$ defined on the submanifold $\Delta$ of $M$. If $\Delta'$ is a face of $\Delta$, then $({\cA^{!!}_{\Delta}})^{!!}_{\Delta'}=\cA^{!!}_{\Delta'}$ (see proposition 2.6 of \cite{mish-oli}) and so we obtain a complex of Lie algebroids on $K$, which will be called the corresponding complex of $\cA$ on the simplicial complex $K$ and denoted by $\{\cA^{!!}_{\Delta}\}_{\Delta\in K}$.

\vspace{3mm}




Let $\underline{\cA}=\{\cA_{\Delta}\}_{\Delta\in K}$ be a complex of Lie algebroids on a simplicial complex $K$. For each simplex $\Delta$ of $K$, the cochain algebra of all smooth forms on $\cA_{\Delta}$ is denoted by $\Omega^{\ast}(\cA_{\Delta};\Delta)$. Let $\Delta$ and $\Delta'$ be two simplices of $K$, with $\Delta'$ face of $\Delta$, and $\varphi_{\Delta,\Delta'}:\Delta'\hookrightarrow \Delta$ the inclusion mapping. By definition of complex of Lie algebroids, we have that $(\cA_{\Delta})^{!!}_{\Delta'}=\cA_{\Delta'}$. Let
$$\varphi^{\cA_{\Delta}}_{\Delta,\Delta'}:\Omega^{\ast}(\cA_{\Delta};\Delta)\longrightarrow \Omega^{\ast}(\cA_{\Delta'};\Delta')$$ be the homomorphism of cochain algebras generated by the inclusion $\varphi_{\Delta,\Delta'}$ (see definition 3.2 of \cite{mish-oli}). Based on Whitney's book \cite{wity-git} and Sullivan's papers \cite{suli-inf} and \cite{suli-Tok}, we give now the definition of piecewise smooth form on a complex of Lie algebroids.


\begin{defn} Let $K$ be a simplicial complex and $\underline{\cA}=\{\cA_{\Delta}\}_{\Delta\in K}$ a complex of Lie algebroids on $K$. A piecewise smooth form of degree $p$ ($p\geq0$) on the complex of Lie algebroids $\underline{\cA}$ is a family $\omega=(\omega_{\Delta})_{\Delta \in K}$ such that, for each $\Delta\in K$, $\omega_{\Delta}\in \Omega^{p}(\cA_{\Delta};\Delta)$ is a smooth
form of degree $p$ on $\cA_{\Delta}$ and, for each $\Delta$, $\Delta'$ $\in K$, with $\Delta'$ face of
$\Delta$, $$\varphi^{\cA_{\Delta}}_{\Delta,\Delta'}(\omega_{\Delta})=\omega_{\Delta'}$$
\end{defn}

\vspace{3mm}

The simplex $\Delta'$ is an embedded compact submanifold of $\Delta$ and so the Lie algebroid $(\cA_{\Delta})^{!!}_{\Delta'}$ can be identified to the Lie algebroid $\im (\varphi_{\Delta,\Delta'})^{!!}$ (proposition 2.4 of \cite{mish-oli}). Hence, for each $x\in \Delta'$, the fibre $(\cA_{\Delta'})_{x}$ is a vector subspace of the fibre $(\cA_{\Delta})_{x}$. The second condition of the definition given above can be stated in the following form: for each $x\in \Delta'$ and vectors $u_{1}$, . . . , $u_{p}$ $\in
(A_{\Delta'})_{x}$ $$\omega_{\Delta'}(x)(u_{1}, \cdot\cdot\cdot, u_{p})=\omega_{\Delta}(x)(u_{1}, \cdot\cdot\cdot, u_{p})$$
Thus, a piecewise smooth form is a collection of smooth forms, each one defined on a transitive Lie algebroid over a simplex of $K$, which are compatible under restrictions to faces. The set of all piecewise smooth forms of degree $p$ on the complex of Lie algebroids $\underline{\cA}$ will be denoted by  $\Omega^{p}(\underline{\cA};K)$. We have then
$$\Omega^{p}(\underline{\cA};K)=\{(\omega_{\Delta})_{\Delta\in
K}:\omega_{\Delta}\in \Omega^{p}(\cA_{\Delta}), \ \ \Delta'\prec\Delta\
\ \Longrightarrow \ \ (\omega_{\Delta})_{/\Delta'}=\omega_{\Delta'}\}$$ in which $(\omega_{\Delta})_{/\Delta'}$ denotes a form $\varphi^{\cA_{\Delta}}_{\Delta,\Delta'}(\omega_{\Delta})$.

Since restrictions of smooth forms are compatible with sums and products, various operations on $\Omega^{p}(\underline{\cA};K)$ can be defined
by the corresponding operations on each simplex of $K$. The set $\Omega^{p}(\underline{\cA};K)$, equipped with these operations, becomes a real vector subspace of $\prod_{\Delta\in K}\Omega^{p}(\cA_{\Delta};\Delta)$, for each natural $p\geq 0$. Thus, $\Omega^{p}(\underline{\cA};K)$ is a module over the subalgebra of $C(|K|;\RR)$ made by all continuous maps $f:|K|\longrightarrow\RR$ that are compatible with restrictions to the faces of $K$ and with smooth restrictions to the faces of $K$. When $p=0$, $\Omega^{0}(\underline{\cA};K)=C(|K|;\RR)$ has a structure of an unitary associative algebra over $\RR$. Moreover, the direct sum $$\Omega^{\ast}(\underline{\cA};K)=\bigoplus_{p\geq 0}\Omega^{p}(\underline{\cA};K)$$ equipped with the exterior product defined by the corresponding exterior product on each algebra $\Omega^{\ast}(\cA_{\Delta};\Delta)=\bigoplus_{p\geq 0}\Omega^{p}(\cA_{\Delta};\Delta)$, is a commutative graded algebra over $\RR$.

In order to obtain a complex of cochains, especially important is the analogues of exterior derivative. This operator also is obtained by the corresponding exterior derivative on each simplex. Namely, if $\underline{\cA}=\{\cA_{\Delta}\}_{\Delta\in K}$ is a complex of Lie algebroids on a
simplicial complex $K$, we can define the mapping $$d^{p}:\Omega^{p}(\underline{\cA};K)\longrightarrow \Omega^{p+1}(\underline{\cA};K)$$
by setting $$d^{p}((\omega_{\Delta})_{\Delta\in K})=(d^{p}_{\Delta}\omega_{\Delta})_{\Delta\in K}$$ for each $\omega=(\omega_{\Delta})_{\Delta\in K}\in \Omega^{p}(\underline{\cA};K)$. Such as in the case of smooth forms on transitive Lie algebroids defined on smooth manifolds (see Kubarski \cite{kuki-dual}), the space $\Omega^{\ast}(\underline{\cA};K)$, with the operations and differentiation above, becomes a commutative differential graded algebra, which is defined over $\RR$.


\begin{defn} \textit{(Piecewise smooth cohomology)}. Keeping the same hypotheses and notations as above, the piecewise smooth cohomology of  $\underline{\cA}$ is the cohomology space of the algebra $\Omega^{\ast}_{\ast}(\underline{\cA};K)$ equipped with the structures defined above. Its cohomology, $H(\Omega^{\ast}(\underline{\cA};K))$, will be denoted by $H^{\ast}(\underline{\cA};K)$.
\end{defn}

\vspace{6mm}

\section{Mayer-Vietoris sequence}

\vspace{3mm}

Assume that $K_{0}$ and $K_{1}$ are two simplicial subcomplexes of a simplicial complex $K$ such that $K=K_{0}\cup K_{1}$. Let $\underline{\cA}=\{\cA_{\Delta}\}_{\Delta\in K}$ be a complex of Lie algebroids on $K$. We can consider the complexes of Lie algebroids $\underline{\cA^{0}}=\{\cA_{\Delta}\}_{\Delta\in K_{0}}$, $\underline{\cA^{1}}=\{\cA_{\Delta}\}_{\Delta\in K_{1}}$ and $\underline{\cA^{0,1}}=\{\cA_{\Delta}\}_{\Delta\in K_{0}\cap K_{1}}$. Our propose now is to know which relations hold
between the cohomology of the complexes of cochains $\Omega^{\ast}(\underline{\cA};K)$, $\Omega^{\ast}(\underline{\cA^{0}};K_{0})$,
$\Omega^{\ast}(\underline{\cA^{1}};K_{1})$ and  $\Omega^{\ast}(\underline{\cA^{0,1}};K_{0}\cap K_{1})$. The answer to this question is the Mayer-Vietoris sequence. We start by defining now the notion of restriction of piecewise smooth forms on a complex of Lie algebroids and state some lemmas on extensions of forms. Once this is done, we establish the result concerned to the Mayer-Vietoris sequence.


\begin{defn} \textit{(Restriction of piecewise smooth forms)} Let $K$ be a simplicial complex and $\underline{\cA}=\{\cA_{\Delta}\}_{\Delta\in K}$ a complex of Lie algebroids on $K$. Let $L$ be a simplicial subcomplex of $K$ and $\underline{\cA}^{L}=\{\cA_{\Delta}\}_{\Delta\in L}$ the
complex of Lie algebroids given by restriction of $\underline{\cA}$ to $L$. If
$\omega=(\omega_{\Delta})_{\Delta\in K}\in\Omega^{p}(\underline{\cA};K)$ is a piecewise smooth form of degree $p$, we can define the
restriction of $\omega$ to the subcomplex $L$, denoted by $\omega_{/L}$, to be the form $$\omega_{/L}=(\omega_{\Delta})_{\Delta\in L}\in\Omega^{p}(\underline{\cA}^{L};L)$$
\end{defn}

\vspace{3mm}

In the conditions of this definition, the equality $$d(\omega_{/L})=(d\omega)_{/L}$$ holds.

\vspace{3mm}

Our first lemma on extensions of piecewise smooth forms is stated for the particular case of a complex of Lie algebroids defined on the simplicial complex made by the canonical simplex and its faces.

\vspace{3mm}

\begin{prop} Let $\Delta_{k}$ denote the canonical $k$-simplex in $\RR^{\infty}$ having the vertices $$e_{0}=(0,0,\dots,0,\dots)$$
$$e_{1}=(1,0,\dots,0,\dots)$$ $$\dots$$ $$e_{k}=(0,0,\dots,1,\dots)$$ ($e_{j}$ is the vector with 1 in the $j$th coordinate and zeros elsewhere).
Let $\cA$ be a transitive Lie algebroid defined on the simplex $\Delta_{k}$ and consider the complexes of Lie algebroids
$\underline{\cA^{\Delta_{k}}}=\{\cA_{\alpha}^{!!}\}_{\alpha\in \Delta_{k}}$ and
$\underline{\cA^{\textrm{\textbf{bd}}\Delta_{k}}}=\{\cA_{\alpha}^{!!}\}_{\alpha\in \textrm{\textbf{bd}}\ \Delta_{k}}$
given by restriction of $\cA$ to the correspondent simplicial complexes $\Delta_{k}$ and $\textrm{\textbf{bd}}\ \Delta_{k}$ respectively. Let
$\xi\in\Omega^{p}\big(\underline{\cA^{\textrm{\textbf{bd}}\Delta_{k}}};\textrm{\textbf{bd}}\ \Delta_{k}\big)$ be a piecewise smooth form
of degree $p$ defined on $\textrm{\textbf{bd}}\  \Delta_{k}$. Then, there is a piecewise smooth form
$\omega\in\Omega^{p}\big(\underline{\cA^{\Delta_{k}}};\Delta_{k}\big)$ of degree $p$ defined on $\Delta_{k}$ such that
$\omega_{/\textrm{\textbf{bd}}\  \Delta_{k}}=\xi$.
\end{prop}

\proofwp. We are going to divide the proof into three parts.


Part 1. This part of the proof follows the proof of the extension lemma 8.3 presented in \cite{gang-diff} and the example (i) and (ii) presented in the seventh section of Sullivan's paper \cite{suli-inf}. Let $\alpha$ be a face of dimension $k-1$ of $\Delta_{k}$, say us, $\alpha$ is the face spanned by the vertices $e_{0}$, $\dots$, $e_{j-1}$, $e_{j+1}$, $\dots$, $e_{k}$. The face $\alpha$ consists of all points $x\in \RR^{\infty}$ such that $$x=t_{0}e_{0}+\dots +t_{j-1}e_{j-1}+0e_{j}+t_{j+1}e_{j+1}+\dots +t_{k}e_{k}$$ with $\sum_{i}t_{i}=1$ and $t_{i}\geq 0$. Let $e_{j}$ be the opposite vertex to the face $\alpha$ and $U$ the complement of this vertex. $U$ is an open subset in $\Delta_{k}$. For each $$x=t_{0}e_{0}+\dots +t_{j-1}e_{j-1}+t_{j}e_{j}+t_{j+1}e_{j+1}+\dots +t_{k}e_{k}\in U$$ we have that $t_{j}\neq 1$ and the element $\frac{t_{0}}{1-t_{j}}e_{0}+\dots+\frac{t_{j-1}}{1-t_{j}}e_{j-1}+\frac{t_{j+1}}{1-t_{j}}e_{j+1}\dots+\frac{t_{k}}{1-t_{j}}e_{k}$ belongs to $\alpha$. So, we may define the map $$\varphi:U\longrightarrow \alpha$$ by
$$\varphi(t_{0}e_{0}+ \dots+t_{j-1}e_{j-1}+t_{j}e_{j}+t_{j+1}e_{j+1}+\dots +t_{k}e_{k})=$$
$$=\frac{t_{0}}{1-t_{j}}e_{0}+\dots+
\frac{t_{j-1}}{1-t_{j}}e_{j-1}+\frac{t_{j+1}}{1-t_{j}}e_{j+1}\dots+
\frac{t_{k}}{1-t_{j}}e_{k}$$ Obviously $\varphi$ is smooth map. Thus,  $\varphi$ is a retraction. Since $\Delta_{k}$ is contractible, the Lie algebroid $\cA$ is a trivial Lie algebroid (see \cite{mish-oli}, proposition 3.5) and then we can find a map $\psi: \cA_{U}\longrightarrow \cA_{\alpha}$ such that $\lambda=(\psi,\varphi)$ is a morphism of Lie algebroids and, for each $x\in \alpha$, $\psi_{x}:\cA_{x}\longrightarrow \cA_{x}$ is the identity map. Consider now a smooth form $\omega\in \Omega^{p}(\cA^{!!}_{\alpha};\alpha)$. Take the form $\lambda^{\ast}\omega$. This form is smooth and belongs to $\Omega^{p}(\cA_{U};U)$. By extension lemmas, the form $\lambda^{\ast}\omega$ damps out smoothly to a smooth form $\widetilde{\omega}\in \Omega^{p}(\cA;\Delta_{k})$. By taking the restriction of $\widetilde{\omega}$ to each face of the simplex $\Delta_{k}$, we obtain a piecewise smooth form $\widetilde{\omega}\in \Omega^{p}(\underline{\cA^{\Delta_{k}}};\Delta_{k}\big)$, which is a piecewise smooth extension of $\omega$.

Part 2. The piecewise smooth form $\widetilde{\omega}$ obtained in the first part has the following property: for each face $\beta$ of $\alpha$,
$$\omega_{/\beta}=0 \ \ \ \Longrightarrow \ \ \ \widetilde{\omega}_{/\beta\ast e_{j}}=0$$ where $\beta\ast e_{j}$ is the join of $\beta$ and the vertex $e_{j}$. This happens because $\widetilde{\omega}_{e_{j}}$ is obviously equal to zero and, for each $x\in \beta\ast e_{j}$ with $x\neq e_{j}$, $\varphi(x)$ lives in $\beta$.

Part 3. Let $\alpha_{0}$, $\dots$, $\alpha_{k}$ be the $k+1$ faces of dimension $k-1$ of $\Delta_{k}$ and let
$$\xi=(\xi_{\alpha})_{\alpha\in \ \textrm{\textbf{bd}}\Delta_{k}}=(\xi_{\alpha_{0}},\xi_{\alpha_{1}},\dots,\xi_{\alpha_{k}})\in \Omega^{p}(\underline{\cA^{\textrm{\textbf{bd}}\Delta_{k}}};\textrm{\textbf{bd}}\ \Delta_{k})$$ be a piecewise smooth form of degree $p$ defined over $\textrm{\textbf{bd}}\Delta_{k}$. By the part 1, the smooth form $\xi_{\alpha_{0}}$ $\in \Omega^{p}(\cA^{!!}_{\alpha_{0}};\alpha_{0})$ can be extended to a smooth form $\widetilde{\xi_{0}}$ $\in \Omega^{p}(\cA;\Delta_{k})$ defined on $\Delta_{k}$ and the form $\widetilde{\xi_{0}}$ defines,
by restriction to each face of $\Delta_{k}$, a piecewise smooth form
$\widehat{\xi_{0}}\in \Omega^{p}\big(\underline{\cA^{\Delta_{k}}};\Delta_{k}\big)$
defined on $\Delta_{k}$. Consider the form $\xi_{1}$ defined by $$\xi_{1}=\xi-\big(\widehat{\xi_{0}}_{\big/\textrm{\textbf{bd}}\Delta_{k}}\big)\in \Omega^{p}\big(\underline{\cA^{\textrm{\textbf{bd}}\Delta_{k}}};\textrm{\textbf{bd}}\ \Delta_{k}\big)$$
The form $\xi_{1}$ vanishes on each point of $\alpha_{0}$. Repeating the same process for the face $\alpha_{1}$ by using the form $\xi_{1}$, the smooth form ${\xi_{1}}_{/\alpha_{1}}\in \Omega^{p}(\cA^{!!}_{\alpha_{1}};\alpha_{1})$ extends to a smooth form $\widetilde{\xi_{1}}\in \Omega^{p}(\cA;\Delta_{k})$ defined on $\Delta_{k}$, and then we obtain a piecewise smooth form $\widehat{\xi_{1}}\in \Omega^{p}\big(\underline{\cA^{\Delta_{k}}};\Delta_{k}\big)$ defined on $\Delta_{k}$ by restrictions to its faces. Since the faces $\alpha_{0}$ and $\alpha_{1}$ have a common vertex, we have that $\widehat{\xi_{1}}_{/\alpha_{0}}=0$ by the second part above. Let
$\xi_{2}=\xi-(\widehat{\xi_{0}}+\widehat{\xi_{1}})_{\big/\textrm{\textbf{bd}}\Delta_{k}}\in \Omega^{p}\big(\underline{\cA^{\textrm{\textbf{bd}}\Delta_{k}}};\textrm{\textbf{bd}}\ \Delta_{k}\big)$.
Then, $\xi_{2}$ is a piecewise smooth form defined on $\textrm{\textbf{bd}}\Delta_{k}$ in which ${\xi_{2}}_{/\alpha_{0}}=0$ and
${\xi_{2}}_{/\alpha_{1}}={\xi_{1}}_{/\alpha_{1}}-\widetilde{\xi_{1}}_{/\alpha_{1}}=0$.
Hence ${\xi_{2}}_{/\alpha_{0}\cup \alpha_{1}}=0$. Therefore, we construct inductively a finite sequence of forms
$\xi_{j}\in\Omega^{p}\big(\underline{\cA^{\textrm{\textbf{bd}}\Delta_{k}}};\textrm{\textbf{bd}}\ \Delta_{k}\big)$ and
$\widehat{\xi_{j}}\in \Omega^{p}_{ps}\big(\underline{\cA^{\Delta_{k}}};\Delta_{k}\big)$, with $j=1,\dots,k+1$, such that
$$\xi_{j}=\xi-(\widehat{\xi_{0}}+\widehat{\xi_{1}}+\dots+\widehat{\xi_{j-1}})_{\big/\partial
\Delta}$$ and $${\xi_{j}}_{/\alpha_{0}\cup \alpha_{1}\cup \dots \cup \alpha_{j-1}}=0$$ for each $j\in \{1,\dots,k+1\}$. Then, setting $$\omega=\widehat{\xi_{0}}+\widehat{\xi_{1}}+\dots+\widehat{\xi_{l}}$$ we have that
$\omega\in \Omega^{p}\big(\underline{\cA^{\Delta_{k}}};\Delta_{k}\big)$ is a piecewise smooth form defined on $\Delta_{k}$ such that
$\omega_{/\textrm{\textbf{bd}}\Delta_{k}}=\xi$ and so the result is proved. {\small $\square$}

\vspace{3mm}

We consider now the extension lemma of piecewise smooth forms on a complex of Lie algebroids defined on a simplicial complex made by any general simplex and its faces.

\vspace{3mm}

\begin{prop} Let $\Delta$ be any simplex of dimension $k$. Let $\cA$ be a transitive Lie algebroid on $\Delta$ and consider the complexes of Lie algebroids $\underline{\cA^{\Delta_{k}}}=\{\cA_{\alpha}^{!!}\}_{\alpha\in \Delta_{k}}$ and
$\underline{\cA^{\textrm{\textbf{bd}}\Delta_{k}}}=\{\cA_{\alpha}^{!!}\}_{\alpha\in \textrm{\textbf{bd}}\ \Delta_{k}}$
given by restriction of $\cA$ to the correspondent simplicial complexes $\Delta_{k}$ and
$\textrm{\textbf{bd}}\ \Delta_{k}$ respectively. Let
$\xi\in\Omega^{p}\big(\underline{\cA^{\textrm{\textbf{bd}}\Delta_{k}}};\textrm{\textbf{bd}}\ \Delta_{k}\big)$
be a piecewise smooth form of degree $p$ defined on $\textrm{\textbf{bd}}\  \Delta_{k}$. Then, there is a piecewise smooth form
$\omega\in\Omega^{p}\big(\underline{\cA^{\Delta_{k}}};\Delta_{k}\big)$ of degree $p$ defined on $\Delta_{k}$ such that
$\omega_{/\textrm{\textbf{bd}}\  \Delta_{k}}=\xi$.
\end{prop}

\vspace{3mm}

\proofwp. There is an affine isomorphism $\varphi$ from the simplex $\Delta_{k}$ onto the simplex $\Delta$ which maps the boundary
$\textrm{\textbf{bd}}\  \Delta_{k}$ onto the boundary $\textrm{\textbf{bd}}\  \Delta$. Consider the transitive Lie algebroid $\varphi^{!!}(\cA)$ on $\Delta_{k}$. Then, $\varphi^{!!}(\cA)$ is isomorphic (non strong isomorphism of Lie algebroids) to the Lie algebroid $\cA$. Take the
inverse image of the form $\xi$, apply the previous proposition, take the direct image and we have the required extension. {\small $\square$}

\vspace{3mm}

In the previous two propositions, we began with a piecewise smooth form defined on whole boundary. However, we can improve slightly last propositions and establish a result concerning extension of piecewise smooth forms when the form is defined, not on all $(k-1)$$-$dimensional faces, but just on some $(k-1)$$-$dimensional faces of $\Delta$. We note this improvement on next proposition.

\vspace{3mm}

\begin{prop} Let $\Delta$ be any simplex of dimension $k$ and $\cA$ a transitive Lie algebroid on the simplex $\Delta$. Consider the complex of Lie algebroids $\underline{\cA^{\Delta}}=\{\cA_{\alpha}^{!!}\}_{\alpha\in \Delta}$ given by restriction of $\cA$ to the correspondent simplicial complex $\Delta$. Suppose that $\alpha_{0}$, $\dots$, $\alpha_{k}$ are the $k+1$ faces of dimension $k-1$ of simplex $\Delta$ and that we have $r$ smooth forms $\xi_{j_{1}}\in \Omega^{p}(\cA^{!!}_{\alpha_{j_{1}}};\alpha_{j_{1}})$, $\dots$, $\xi_{j_{r}}\in \Omega^{p}(\cA^{!!}_{\alpha_{j_{r}}};\alpha_{j_{r}})$ with $\{j_{1},\dots,j_{r}\}\subset \{1,\dots,k\}$ such that, for each $j_{i}$, $j_{e}$ with $\alpha_{j_{i}}\cap \alpha_{j_{e}}$ non empty, the forms $\xi_{j_{i}}$ and $\xi_{j_{e}}$ agree on the intersection $\alpha_{j_{i}}\cap \alpha_{j_{e}}$. Then, there is a piecewise smooth form $\omega\in\Omega^{p}_{ps}\big(\underline{\cA^{\Delta}};\Delta\big)$ such that $\omega_{/\alpha_{j_{i}}}=\xi_{j_{i}}$ for $i=1,\dots,r$.
\end{prop}

\vspace{3mm}

\proofwp. For each vertex which does not belong to any face $\alpha_{j_{1}}$, $\dots$, $\alpha_{j_{r}}$, we take the smooth form zero on this vertex and then we have a family of smooth forms, each one defined on each vertex of $\Delta$. Now, for each two vertices defining a face of $\Delta$ of dimension 1 different of any $1$$-$dimensional face of $\alpha_{j_{i}}$, with $i=1,\dots,r$, we apply the previous proposition and we obtain a piecewise smooth form defined on the skeleton of $\Delta$ of dimension 1. This piecewise smooth form is an extension of each smooth forms given on the $1$$-$dimensional faces of $\alpha_{j_{i}}$ ($i=1,\dots,r$) by restriction of the forms $\xi_{j_{r}}$. We repeat the same argument for dimension 2. This process will end on dimension $k$ and the form obtained is a piecewise smooth form defined on $\Delta$, which is an extension of the forms $\xi_{j_{i}}$ ($i=1,\dots,r$). {\small $\square$}

\vspace{3mm}

From last propositions, we easily obtain the next general lemma on extensions of piecewise smooth forms on complexes of Lie algebroids.

\vspace{3mm}

\begin{prop} Let $K$ be a simplicial complex and $\underline{\cA}=\{\cA_{\alpha}\}_{\alpha\in K}$ a complex of Lie algebroids
on $K$. Let $L$ be a simplicial subcomplex of $K$ and consider the subcomplex of Lie algebroids $\{\cA_{\alpha}\}_{\alpha\in L}$ defined
on $L$. Then, any piecewise smooth form of degree $p$ defined on $L$ can be piecewise smoothly extended to a piecewise smooth form of
degree $p$ defined on the whole $K$.
\end{prop}

\vspace{3mm}

Next proposition concerns the short exact sequence which generates the Mayer-Vietoris sequence in cohomology. For transitive Lie algebroids on smooth manifolds, the short exact sequence is presented in the third section of \cite{kuki-dual}.

\vspace{3mm}

\begin{prop} Let $K$ be a simplicial complex and $\underline{\cA}=\{\cA_{\alpha}\}_{\alpha\in K}$ a complex of Lie algebroids on $K$.
Let $K_{0}$ and $K_{1}$ be two simplicial subcomplexes of $K$ such that $K=K_{0}\cup K_{1}$ and set $L=K_{0}\cap K_{1}$. Consider the complexes of Lie algebroids $\underline{\cA^{0}}=\{\cA_{\alpha}\}_{\alpha\in K_{0}}$, $\underline{\cA^{1}}=\{\cA_{\alpha}\}_{\alpha\in K_{1}}$ and
$\underline{\cA^{0,1}}=\{\cA_{\alpha}\}_{\alpha\in L}$ given by restriction of $\cA$ to the simplicial subcomplexes $K_{0}$, $K_{1}$ and $L$. Then, it holds a exact short sequence of cochain complexes
$$\{0\}\longrightarrow \Omega^{\ast}(\underline{\cA};K)\mapr{\lambda^{\ast}}\Omega^{\ast}(\underline{\cA^{0}};K_{0})\oplus
\Omega^{\ast}(\underline{\cA^{1}};K_{1})
\mapr{\mu^{\ast}}\Omega^{\ast}(\underline{\cA^{0,1}};L)\longrightarrow
\{0\}$$

\noindent in which the linear maps $$\lambda^{p}:
\Omega^{\ast}(\underline{\cA};K)\longrightarrow
\Omega^{\ast}(\underline{\cA^{0}};K_{0})\oplus
\Omega^{\ast}(\underline{\cA^{1}};K_{1})$$
$$\mu^{p}:\Omega^{\ast}(\underline{\cA^{0}};K_{0})\oplus
\Omega^{\ast}(\underline{\cA^{1}};K_{1})\longrightarrow
\Omega^{\ast}(\underline{\cA^{0,1}};L)$$

\noindent are defined by
$\lambda^{p}(\omega)=(\omega_{/K_{0}},\omega_{/K_{1}})$ and
$\mu^{p}(\xi,\eta)=\eta_{/L}-\xi_{/L}$.
\end{prop}

\vspace{3mm}

\proofwp. As in the case of smooth forms on a transitive Lie algebroid over a smooth manifold, the exterior derivative commutes with the restrictions to a simplicial subcomplexes and, since $d^{p}(\xi,\eta)=(d^{p}(\xi),d^{p}(\eta))$, one deduces immediately that $\lambda^{\ast}$ and $\mu^{\ast}$ are effectively cochain complex morphisms. Obviously, the linear map $\lambda^{p}$ is injective. Since, for each piecewise smooth form $\omega\in
\Omega^{p}(\underline{\cA};K)$, the forms $\omega_{/K_{0}}$ and $\omega_{/K_{1}}$ have the same restriction $\omega_{L}$ to $L$, we conclude that $\mu^{p}\circ \lambda^{p}=0$, and hence the image of the linear map $\lambda^{p}$ is contained in the kernel of the linear map $\mu^{p}$. Reciprocally, if $\mu^{p}(\xi,\eta)=0$, we have $\xi_{\alpha}=\eta_{\alpha}$, for each $\alpha\in L$, and this equality allows to define a piecewise smooth form $\omega\in \Omega^{p}(\underline{\cA};K)$ by the condition $\omega_{\alpha}=\xi_{\alpha}$, for each $\alpha\in K_{0}$, and $\omega_{\alpha}=\eta_{\alpha}$, for each $\alpha\in K_{1}$. We have then $\lambda^{p}(\omega)=\mu^{p}(\xi,\eta)$. We want now to prove that $\mu^{p}$ is surjective. Let $\gamma \in \Omega^{p}(\underline{\cA^{0,1}};L)$ be a piecewise smooth form and consider the piecewise smooth form $-\frac{1}{2}\gamma \in \Omega^{p}(\underline{\cA^{0,1}};L)$. By the extension lemma, we can consider a piecewise smooth form $\alpha \in
\Omega^{p}(\underline{\cA^{0}};K_{0})$ such that $\alpha_{/L}=-\frac{1}{2}\gamma$. Analogously, we can consider a piecewise smooth form $\beta \in \Omega^{p}(\underline{\cA^{1}};K_{1})$ such that $\beta_{/L}=\frac{1}{2}\gamma$. We have then that $\mu^{p}(\alpha,\beta)=\gamma$. {\small $\square$}

\vspace{3mm}

By applying the zig-zag lemma to the sequence above, we obtain the long exact sequence in cohomology

\vspace{3mm}

$$\begin{CD}
\cdots @>>> H^{p-1}(\underline{\cA^{0,1}};L) @>\partial^{p-1}>> H^{p}(\underline{\cA};K) @>H^{p}(\lambda^{\ast})>> H^{p}(\underline{\cA^{0}};K_{0})\oplus H^{p}(\underline{\cA^{1}};K_{1})
\end{CD}$$

$$\begin{CD}
H^{p}(\underline{\cA^{0}};K_{0})\oplus H^{p}(\underline{\cA^{1}};K_{1}) @>H^{p}(\mu^{\ast})>> H^{p}(\underline{\cA^{0,1}};L) @>\partial^{p}>> H^{p+1}(\underline{\cA};K)  @>>> \cdots \\
\end{CD}$$

\vspace{3mm}

\noindent which is the Mayer-Vietoris sequence for piecewise smooth cohomology of complexes of Lie algebroids.

\vspace{3mm}


\end{document}